\documentstyle[amscd,amssymb,verbatim]{amsart}
\newcommand{\vareps}{\varepsilon}

\renewcommand{\mod}{\operatorname{mod}}

\newcommand{\Tr}{\operatorname{Tr}}

\newcommand{\OO}{{\cal O}}

\newcommand{\DD}{{\cal D}}

\newcommand{\Gal}{\operatorname{Gal}}

\newcommand{\SL}{\operatorname{SL}}
\newcommand{\G}{{\Bbb G}}

\newcommand{\HHom}{{\cal H}om}

\newcommand{\Coh}{\operatorname{Coh}}

\newcommand{\Spec}{\operatorname{Spec}}

\renewcommand{\Sp}{\operatorname{Sp}}

\renewcommand{\P}{{\Bbb P}}

\newcommand{\si}{\sigma}

\newcommand{\Pic}{\operatorname{Pic}}

\newcommand{\de}{\delta}
\newcommand{\eps}{\epsilon}

\renewcommand{\ker}{\operatorname{ker}}

\numberwithin{equation}{section}

\newtheorem{thm}{Theorem}[section]
\newtheorem{prop}[thm]{Proposition}
\newtheorem{lem}[thm]{Lemma}
\newtheorem{cor}[thm]{Corollary}
\newenvironment{rem}{\vspace{3mm}\noindent
{\bf Remark.}}{\vspace{3mm}}
\newenvironment{defi}{\vspace{3mm}\noindent
{\bf Definition.}}{\vspace{3mm}}
\newenvironment{rems}{\vspace{3mm}
\noindent {\bf Remarks.}}{\vspace{3mm}}

\newcommand{\Pf}{\noindent {\it Proof}}

\newcommand{\ov}{\overline}

\newcommand{\rk}{\operatorname{rk}}

\newcommand{\MM}{{\cal M}}
\newcommand{\TT}{{\cal T}}

\newcommand{\UU}{{\frak U}}

\newcommand{\Hom}{\operatorname{Hom}}

\renewcommand{\a}{\alpha}

\newcommand{\om}{\omega}
\newcommand{\De}{\Delta}
\newcommand{\la}{\lambda}

\newcommand{\C}{{\Bbb C}}

\newcommand{\Z}{{\Bbb Z}}
\newcommand{\Q}{{\Bbb Q}}
\newcommand{\La}{\Lambda}

\newcommand{\wt}{\widetilde}

\newcommand{\sub}{\subset}
\newcommand{\ed}{\qed\vspace{3mm}}

\newcommand{\GL}{\operatorname{GL}}
\newcommand{\tr}{\operatorname{tr}}
\newcommand{\Spin}{\operatorname{Spin}}
\newcommand{\SO}{\operatorname{SO}}

\title{$K$-theoretic exceptional collections at roots of unity}
\author{A. Polishchuk}
\address{Department of Mathematics, University of Oregon, Eugene, OR 97405}
\email{apolish@@math.uoregon.edu}
\thanks{Supported in part by NSF grant}

\begin{document}
\begin{abstract} Using cyclotomic specializations of the
equivariant $K$-theory with respect to a torus action
we derive congruences for discrete invariants of exceptional objects in derived categories
of coherent sheaves on a class of varieties that includes Grassmannians and smooth quadrics.
For example, we prove that if $X=\P^{n_1-1}\times\ldots\times\P^{n_k-1}$, where
$n_i$'s are powers of a fixed prime number $p$, then the rank of an exceptional object 
on $X$ is congruent to $\pm 1$ modulo $p$. 
\end{abstract}
\maketitle

\bigskip

\section{Introduction}

This paper is concerned with $K$-theory classes of exceptional objects in the derived category
of coherent sheaves $D(X):=D^b(\Coh X)$ on a smooth projective variety $X$. 
Recall that an object $E$ of a $k$-linear triangulated category is called {\it exceptional} if
$R\Hom(E,E)=k$. An {\it exceptional collection} is a collection of exceptional objects
$(E_1,\ldots,E_n)$ such that $R\Hom(E_i,E_j)=0$ for $i>j$. An exceptional collection is
called {\it full} if it generates the entire triangulated category. For example, 
$(\OO,\OO(1),\ldots,\OO(n))$ is a full exceptional collection in $D(\P^n)$ (see \cite{Be}). 
We refer to \cite{GK} for more background on exceptional collections (see also section 3.1 of 
\cite{Bridge} for a brief introduction). 
There is a naturally defined action of the braid group $B_n$ on the set of exceptional collections of length $n$ given by mutations. 
It is conjectured that in the case of full exceptional collections this action is transitive (where each object is considered up to a shift $E\mapsto E[m]$). We refer to this property
as {\it constructibility}. This property is known only in some low-dimensional cases (e.g., it is checked
for Del Pezzo surfaces in \cite{KO}).  Note that if $X$ is a smooth projective variety then
for a full exceptional collection $(E_1,\ldots,E_n)$ in $D(X)$ 
the classes $([E_i])$ in the Grothendieck group
$K_0(X)$ form a basis over $\Z$.
Furthermore, this basis is semiorthogonal with respect to the Euler bilinear form 
$$\chi([V],[W]):=\chi(X,V^*\otimes W)$$
on $K_0(X)$, i.e., we have $\chi([E_i],[E_j])=0$ for $i>j$, $\chi([E_i],[E_i])=1$. One still has
an action of the braid group on the set of semiorthogonal bases in $K_0(X)$, so the problem
of constructibility can be formulated at the level of $K_0(X)$. Even for this question very little is
known (\cite{Nogin} seems to be the only work dealing with $3$-dimensional examples). 
For example, this problem is open for projective spaces of dimension $\ge 4$.

In the case when $X$ admits an action by an algebraic torus $T=\G_m^n$,
every exceptional object in $D(X)$ can be equipped with a $T$-equivariant structure, 
i.e., comes from an object of $D^b(\Coh^T(X))$ (see Lemma \ref{equiv-lem}).
Thus, the constructibility question
can also be asked for bases in the $T$-equivariant $K$-group $K_0^T(X)$, semiorthogonal
with respect to the equivariant Euler form. 
Note that $K_0^T(X)$ is a module over the representation ring 
$R=R(T)\simeq\Z[x_1^{\pm 1},\ldots,x_n^{\pm 1}]$. The main observation of this paper is that in some situations one can choose carefully an element $t_0\in T$ of finite order $N$ such that
after the specialization with respect to the homomorphism $\Tr(t_0,?):R\to\Z[\sqrt[N]{1}]$
the equivariant Euler form becomes Hermitian (and positive-definite). 
This means that every full exceptional collection 
provides an orthonormal basis of $K_0^T(X)\otimes_{R}\Z[\sqrt[N]{1}]$ with respect to the
specialization of the Euler form. The action of the braid group in this specialization reduces to
the action by permutations of basis vectors (up to a sign). Hence, the analog of constructibility
at this level should assert that the obtained orthonormal basis does not depend on a full
exceptional collection up to rescaling\footnote{Rescaling by a root of unity corresponds to changing a $T$-equivariant structure on an exceptional object. An additional sign may appear as en effect of
mutations.}. We observe that this is indeed the case because of
the following simple consequence of Kronecker's theorem (see Proposition \ref{number-prop}):
if $z_1,\ldots,z_n$ is a set of cyclotomic integers such that $|z_1|^2+\ldots+|z_n|^2=1$
then all $z_i$'s are zero except one which is a root of unity.


Examples of the above situation include some homogeneous spaces (e.g., Grassmannians and quadrics), Hirzebruch surfaces $F_n$ for even $n$, as well as products of such varieties. In the case when $N$ is a power 
of a prime $p$ we can make further specialization to $K_0(X)\otimes\Z/p\Z$. This way we get
congruences modulo $p$ for classes of exceptional objects (see Corollary \ref{rk-cor} and Theorem
\ref{congr-proj-thm} for the case of projective spaces and their products).

Here is a precise statement of our main result.

\begin{thm}\label{main-thm} 
Let $X$ be a smooth projective variety over $\C$ equipped with an action of an algebraic torus
$T$. Assume 
that the set of invariant points $X^T$ is finite and there exists an element
of finite order $t_0\in T$ such that 

\noindent
$(\star)$ for every $p\in X^T$ one has $\det(1-t_0, T_pX)\neq 0$ 
and $\det(t_0, T_pX)=(-1)^{\dim X}$, where $T_pX$ is the tangent space at $p$. 

\noindent
Consider the homomorphism
$$\Tr(t_0,?): R\to\Z[\sqrt[N]{1}],$$
where $N$ is the order of $t_0$, and set $K_{t_0}=K^T_0(X)\otimes_R \Z[\sqrt[N]{1}]$.
Then

\vspace{2mm}

\noindent
(i) the equivariant Euler form $\chi^T(\cdot,\cdot)$ on $K^T_0(X)$ specializes to a Hermitian form
on $K_{t_0}$ with values in $\Z[\sqrt[N]{1}]$.

\vspace{2mm}

\noindent
(ii) If $E$ is an exceptional object of $D(X)$ equipped with a
$T$-equivariant structure (i.e., a lift to an object of $D^b(\Coh^T X)$)
then the class $v(E)$ of $E$ in $K_{t_0}$
has length $1$ with respect to the Hermitian form in (i). 
If $(E_1,E_2)$ is an exceptional pair in $D(X)$, where both
$E_1$ and $E_2$ are equipped with a $T$-equivariant structure,
then the vectors $v(E_1)$ and $v(E_2)$ are orthogonal with respect to this form.

\vspace{2mm}

\noindent
(iii) Assume $D(X)$ admits a full exceptional collection $(E_1,\ldots,E_n)$ 
where each $E_i$'s is equipped with a $T$-equivariant structure, and let
$(v_1=v(E_1), \ldots, v_n=v(E_n))$ be the corresponding orthonormal $\Z[\sqrt[N]{1}]$-basis
of $K_{t_0}$. Then every unit vector in $K_{t_0}$
is of the form $\pm\zeta v_i$ for some $i$ and some $N$th root of unity $\zeta$.

\vspace{2mm}

\noindent
(iv) In the situation of (iii) 
assume in addition that the action of $T$ on $X$ extends to an action of an algebraic group 
$N\supset T$, such that $T$ is a normal subgroup in $N$ and for some element
$w\in N$ one has $w t_0 w^{-1}=t_0^{-1}$. 
Assume in addition that all the objects $E_i$ admit a $N$-equivariant
structure. Then for every exceptional object $E$ of $D(X)$ admitting an $N$-equivariant
structure one has $v(E)=\pm v_i$ for some $i$. 
\end{thm}

\begin{rems}
1. We will determine for which generalized Grassmannians $G/P$ (where $P$ is a maximal
parabolic subgroup) there exists an element $t_0$ satisfying $(\star)$ in Proposition 
\ref{cases-prop}, leaving out only several cases with $G$ of type $E_7$ and $E_8$. 
In particular, for $G$ of classical type the spaces we get are either Grassmannians, or
smooth quadrics or (connected components) of maximal isotropic Grassmannians, orthogonal
or symplectic.

\noindent
2. The assumptions of part (iv) are often easy to check when $T$ is a maximal torus of
a simply connected semisimple group $G$ acting on $X$: it suffices to find an element of the Weyl group $W$ sending $t_0$ to $t_0^{-1}$ (since in this case every exceptional object admits a $G$-equivariant structure by Lemma \ref{equiv-lem}). In almost all of the cases considered in Proposition \ref{cases-prop} 
this holds for $w_0$, the element of maximal length in $W$ (the exception is the case of type $A_n$
where one should take a different element of $W$).

\noindent
3. Assume that $(X',T',t'_0\in T')$ is another data such that the assumption $(\star)$
of Theorem \ref{main-thm} is satisfied.
Let equip $X\times X'$ with the natural action of the torus $T\times T'$. Then the element
$(t_0,t'_0)\in T\times T'$ will still satisfy the assumption $(\star)$.
\end{rems}

In the case when the order of $t_0$ in the above Theorem is a power of prime we can
derive some congruences in the usual Grothendieck group of $X$.

\begin{cor}\label{rk-cor} 
In the situation of Theorem \ref{main-thm}(iii) assume in addition that the order of $t_0$
equals $p^k$, where $p$ is a prime. Then the reduction of the Euler form $\chi(\cdot,\cdot)$
modulo $p$ is symmetric. Furthermore,
for every exceptional object $E$ of $D(X)$ one has the following congruence in $K_0(X)\otimes\Z/p\Z$:
\begin{equation}\label{mod-p-congr}
[E]\equiv \pm[E_i]
\end{equation} 
for some $i\in [1,n]$.
Also, for such an object one has
\begin{equation}\label{sum-chi-congr}
\sum_{i=1}^n\chi(E_i,E)\equiv \pm 1\mod(p).
\end{equation}
\end{cor}

\begin{rem}
Note that in the situation of the above Corollary there are typically more
vectors of length $1$ in $K_0(X)\otimes\Z/p\Z$ than just those coming from exceptional objects, 
so the congruence \eqref{mod-p-congr} cannot be obtained by just looking at 
$K_0(X)\otimes\Z/p\Z$.
On the other hand, in the case $p=2$ the fact that the ranks of exceptional objects are odd 
can often be checked only with the help of $K_0(X)\otimes\Z/2\Z$ --- see
Remark after Corollary \ref{odd-prod-cor}.
\end{rem}

In most of our examples the torus $T$ is a maximal torus in a connected reductive group $G$ acting
on $X$ in such a way that the center $Z_G\sub G$ acts trivially. In this situation
we have a decomposition of the category $\Coh^G(X)$
of $G$-equivariant coherent sheaves (and of its derived category) into the direct sum of subcategories
indexed by characters of $Z_G$. We say that an object $V\in D(X/G):=D^b(\Coh^G(X))$ is  
{\it central} if it belongs to one of these subcategories, and we call the coresponding character
$\chi:Z_G\to\G_m$ the central character of $V$ (when $V$ is a $G$-equivariant coherent sheaf
this means that $Z_G$ acts on $V$ through $\chi$). 
For example, any indecomposable object in $D(X/G)$ is central. 
Note that in the case when the commutator of $G$ is simply connected every exceptional object
in $D(X)$ admits a $G$-equivariant structure (see Lemma \ref{equiv-lem}), 
and hence can be viewed as a central object of $D(X/G)$.
Using a $G$-equivariant structure often allows to extract more precise information on the
class of an exceptional object in $K_0^T(X)$ (see Propositions \ref{central-orth-prop},
\ref{Galois-prop}).


In section \ref{appl-sec} we will 
consider some concrete examples when the conditions of Theorem \ref{main-thm}
are satisfied. In the case of projective spaces we will prove the following result.

\begin{thm}\label{congr-proj-thm}
(i) Let $p$ be a prime, $n=p^k$, and let 
$V$ be a central object in $D(\P^{n-1}/\GL_{n})$.
Then one has the following congruence in $K_0(\P^{n-1})\otimes \Z/p\Z$:
$$[V]\equiv \cases 0, & \rk(V)\equiv 0\mod(p),\\
\rk(V)[\OO(\deg(V)/\rk(V))], & \rk(V)\not\equiv 0\mod(p),\endcases$$
where we use the fact that the class of $\OO(m)$ in 
$K_0(\P^{n-1})\otimes\Z/p\Z$ depends only on $m\mod(n)$.

\noindent
(ii) If in the above situation $E$ is an exceptional object of $D(\P^{n-1})$ then
$\rk(E)\equiv\pm 1 \mod(p)$.
If $(E_1,E_2)$ is an exceptional pair then 
$$\frac{\deg(E_1)}{\rk(E_1)}\not\equiv \frac{\deg(E_2)}{\rk(E_2)} \mod(p^k),$$
and $\chi(E_1,E_2)\equiv 0 \mod(p)$. 
For a full exceptional collection $(E_1,\ldots,E_{n})$ in $D(\P^{n-1})$ the slopes modulo $p^k$
$(\deg(E_i)/\rk(E_i) \mod (p^k))$ form a complete system of remainders modulo $p^k$.

\noindent
(iii) Let $p$ be a prime and let $X=\P^{n_1-1}\times\ldots\times\P^{n_k-1}$, where $n_i=p^{k_i}$.
Then for an exceptional object $E\in D(X)$ the class of $E$ in 
$K_0(X)\otimes\Z/p\Z$ coincides up to a sign with
the class of some line bundle. In particular, $\rk(E)\equiv\pm 1\mod(p)$.
\end{thm}

There is a relative version of this result for the product of a projective space with a smooth
projective variety (see Theorem \ref{rel-proj-thm}).
We will also consider other situations where Theorem \ref{main-thm}
can be applied. For example, we will prove that if $p$ is a prime then the rank of an exceptional
object on the product of Grassmannians $G(k_1,p)\times\ldots\times G(k_s,p)$ is not divisible by $p$.
(see Theorem \ref{Grass-thm}). In the case of Grassmannians and smooth quadrics
we derive some relations between the rank and the central character of an exceptional object.
For example, we prove that an exceptional object on the smooth quadric of dimension
$2^r$ (with $r\ge 2$)
can be equipped with an $\SO(2^r+2)$-equivariant structure if and only if it has an odd rank
(see Proposition \ref{even-quadric-prop}). For maximal isotropic Grassmannians (orthogonal and 
symplectic) full exceptional collections have been constructed only in few cases, so our results are more
limited. For example, we show that all exceptional objects in the derived category of $OG(5,10)$ (resp., 
$SG(3,6)$) have an odd rank (see Theorem \ref{isotropic-thm}). As examples of non-homogeneous varieties
we consider Hirzebruch surfaces $F_n$ with even $n$.
We show that Theorem \ref{main-thm} applies in this case with $N=4$, and hence
the class of every exceptional object in $K_0(F_n)\otimes\Z/2\Z$ coincides 
with one of the $4$ classes coming from line bundles. 
Similar result holds for products of such surfaces, as well as for their products with
projective spaces of dimension $2^k-1$ (see Corollary \ref{odd-prod-cor}).
In particular, the rank of every exceptional object on such products is odd. 

\bigskip

\noindent
{\it Notations and conventions}.
We work with smooth projective varieties over $\C$. 
We denote by $\zeta_n$ a primitive $n$-th root of unity in $\C$ and by $\Z[\sqrt[n]{1}]\sub \C$
the subring generated by $\zeta_n$.
$R(G)$ denotes the representation ring (over $\Z$) of an algebraic group $G$.
For an algebraic group $G$ acting on a variety $X$ we always assume the existence
of a $G$-equivariant ample line bundle on $X$. 
We denote by $\Coh^G(X)$ (resp., $\Coh(X)$) the category of $G$-equivariant (resp., usual)
coherent sheaves on $X$ and
by $D(X/G)$ (resp., $D(X)$) its bounded derived category.

\section{Lefschetz formula in equivariant $K$-theory and the proof of Theorem \ref{main-thm}}

Let $X$ be a smooth projective variety equipped with an action of an algebraic group $G$.
We denote by $K_0^G(X)$ the Grothendieck group of the category of $G$-equivariant
coherent sheaves on $X$. We always assume that $X$ admits a $G$-equivariant ample line bundle
(this is automatic if $G$ is linear algebraic, see \cite{Kamb}).
In this case every $G$-equivariant coherent sheaf admits a finite resolution by
$G$-equivariant bundles, so $K_0^G(X)$ can also be defined using $G$-equivariant bundles. 
We can view the group $K_0^G(X)$ as a module over $R(G)$ in a natural way. It is also
equipped with the commutative product induced by the tensor product and with
the involution $[V]\mapsto [V^*]$, where $V^*$ is the dual vector bundle to $V$.
Since cohomology of a $G$-equivariant coherent sheaf are equipped with $G$-action, 
we obtain the $G$-equivariant Euler characteristic
$$\chi^G(V)=\chi^G(X,V)=\sum_{i\ge 0}(-1)^i [H^i(X,V)]$$
with values in $R(G)$. Similarly, we have an equivariant version of
the Euler bilinear form with values in $R(G)$ defined by
$$\chi^G(V,W)=\chi^G([V^*\otimes W]).$$
An object $V$ of $D(X)=D^b(\Coh X)$ (resp., of $D(X/G)=D^b(\Coh^G X)$) has an associated
class $[V]$ in $K_0(X)$ (resp., $K_0^G(X)$). By a {\it $G$-equivariant structure} on an object
$V\in D(X)$ we mean an object $\wt{V}\in D(X/G)$ such that the image of $\wt{V}$ under
the forgetful functor $D(X/G)\to D(X)$ is isomorphic to $V$. 

\begin{lem}\label{exc-basis-lem} 
Let $(E_1,\ldots,E_n)$ be a full exceptional collection in $D(X)$.
Assume that each $E_i$ is equipped with a $G$-equivariant structure. Then the classes
$([E_i])$ form a basis of $K^G_0(X)$ as $R(G)$-module.
\end{lem}

\Pf . Recall that with every exceptional object $E$ in $D(X)$ one associates a functor
$L_E$ from $D(X)$ to itself, so that there is a distinguished triangle (functorial in $F$)
\begin{equation}\label{LEF-eq}
L_E(F)\to R\Hom(E,F)\otimes E\to F\to L_E(F)[1].
\end{equation}
Moreover, we have the distinguished triangle of the corresponding kernels in $D(X\times X)$:
$$K\to E^*\boxtimes E\to \De_*\OO_X\to K[1],$$
where $E^*=R\HHom(E,\OO_X)$.
Now assume that $E$ is equipped with a $G$-equivariant structure.
Then we can represent $E$ by a complex of $G$-equivariant vector bundles.
Note that the canonical morphism $E^*\boxtimes E\to \De_*\OO_X$ of complexes of sheaves on 
$X\times X$ is compatible with the diagonal action of $G$.
Hence, there is an analog of \eqref{LEF-eq} in $D(X/G)$ with
$R\Hom(E,F)$ replaced by $R\pi_*(E^*\otimes F)\in D(R(G)-\mod)$, where $\pi:X\to \Spec(\C)$
is the structure morphism. The operator on $K_0^G(X)$ corresponding to the functor $F\mapsto L_E(F)$ 
is given by $[F]\mapsto \chi^G([E],[F])\cdot [E]-[F]$. Hence, the fact that $L_{E_1}\ldots L_{E_n}F=0$
for every $F$ translates into the assertion that the classes $([E_i])$ span $K_0^G(X)$ over
$R(G)$. Using the semiorthogonality condition $\chi^G(E_i,E_j)=0$ for $i>j$, $\chi^G(E_i,E_i)=1$
for all $i$, one easily checks linear independence of these classes over $R(G)$.
\ed

Thus, it is important to know when exceptional objects can be equipped with equivariant
structures. The following result shows that this is always possible for connected reductive groups $G$
with simply connected commutant.

\begin{lem}\label{equiv-lem}
(i) Let $X$ be  a smooth projective variety equipped with an action of a linear algebraic group $G$, and let
$E$ be an exceptional object of $D(X)$. 
Assume that $G$ has trivial Picard group and has no nontrivial central
extensions by $\G_m$ in the category of algebraic groups.
Then $E$ admits a $G$-equivariant structure, 
unique up to tensoring with a character of $G$.

\noindent
(ii) The above assumptions on $G$ are satisfied if $G$ is a connected 
reductive group with $\pi_1(G)$ torsion free, or equivalently, with simply connected
commutant $\DD G$.
\end{lem}

\Pf . (i) We will use the algebraic stack $\MM$ ``parametrizing" objects $E\in D(X)$ with
$\Hom^i(E,E)=0$ for $i<0$ (see \cite{Lieblich}, we could also use the stack defined in \cite{Inaba}).
More precisely, in the terminology of \cite{Lieblich} $\MM$ represents the functor of universally
glueable relatively perfect complexes. 
Consider the pull-back $a^*E$ of $E$ via the action map $a:G\times X\to X$.
By Proposition 2.19 of \cite{Lieblich}, $a^*E$ is a family over $G$ of the above type.
Hence, we have the corresponding morphism $G\to\MM$.
Since the tangent space to $\MM$ at $E$ is $\Hom^1(E,E)=0$, it follows
that this morphism is constant, so the objects
$a^*E$ and $p_2^*E$ on $G\times X$ become isomorphic over a covering of $G$ in flat topology.
Note that if $U\to G$ is one of the elements of this cover then automorphisms
of $p_2^*E$ over $U$ reduce to mutliplication by invertible functions on $U$.
Now the triviality of $\Pic(G)$ implies that we can choose a global isomorphism
$\a:a^*E\to p_2^*E$. Next, we should try to check the cocycle condition for $\a$
on $G\times G\times X$. The obstacle will be some group $2$-cocycle of $G$ with values in $\G_m$.
By the assumptions $G$ has no central extensions by $\G_m$, hence we can adjust $\a$ so that
the cocycle condition will be satisfied. It remains to use the argument of Theorem 3.2.4 of \cite{BBD} to deduce that $E$ can be represented by a complex of $G$-equivariant sheaves. More precisely, 
we observe that $G$-equivariant sheaves of $\OO$-modules on $X$ can be viewed as strict simplicial
systems of sheaves of $\OO$-modules over the simplicial system $(G^n\times X)$ associated with
the action of $G$ on $X$. The above construction extends $E$ to a strict simplicial system in
the corresponding derived category. Now we can apply the argument 
of Theorem 3.2.4 of \cite{BBD} observing that the vanishing condition used in Proposition 2.9 of
{\it loc. cit.} boils down to the vanishing of $\Hom^i(E,E)$ for $i<0$ (similar vanishing for the pull-back
of $E$ to $G^n\times E$ follows by the K\"unneth formula). 

The uniqueness part is checked as follows.
Suppose $E$ and $E'$ are objects of $D(X/G)$ that become isomorphic in $D(X)$ and assume that 
$\Hom_{D(X)}(E,E)=\C$. Then $\Hom_{D(X)}(E,E')$ is a one-dimensional representation of $G$. Hence,
after tensoring $E'$ with a character of $G$ we will get a morphism $f:E\to E'$ in $D(X/G)$ such that
$f$ induces an isomorphism in $D(X)$, so $f$ itself is an isomorphism. 

\noindent
(ii) The equivalence of two conditions on $G$ follows from Corollary 1.7 of \cite{Mer}.
The triviality of the Picard group follows from Proposition 1.10 of \cite{Mer} (see also \cite{Popov}). 
Now let $1\to \G_m\to \wt{G}\to G\to 1$ be a central extension of
$G$ by $\G_m$. Then the derived
group $\DD\wt{G}$ of $\wt{G}$ is a connected semisimple group and we have an isogeny 
$\DD\wt{G}\to \DD G$. Since $\DD G$ is simply connected this implies that the above extension
splits over $\DD G\sub G$. Therefore, it is induced by a central extension of the torus
$G/\DD G$ by $\G_m$.
Hence, we are reduced to the case when $G$ is a torus. 
In this case $\wt{G}$ is a connected solvable algebrac
group with trivial unipotent radical, hence, $\wt{G}$ is itself a torus. Therefore, the above sequence
splits. 
\ed

Let $X$ be a smooth projective variety equipped with an action of an algebraic torus $T\simeq\G_m^n$ 
that has a finite number of stable points $X^T$, and let $R=R(T)$.
The usual (non-equivariant) $K$-group can be recovered from the equivariant one due to 
the following result of Merkurjev.

\begin{lem}\label{eval-lem}(\cite{Mer}, Cor.4.4) 
The natural map
$$K_0^T(X)\otimes_{R}\Z\to K_0(X),$$
induced by the homomorphism $\Tr(1,?):R\to\Z$,
is an isomorphism.
\end{lem}

Consider the natural map of $R$-modules
\begin{equation}\label{loc-map1}
K^T_0(X)\to\bigoplus_{p\in X^T}R
\end{equation}
given by the restriction to $X^T$. It is well known that it
becomes an isomorphism after tensoring over $R$ with the quotient field of $R$
(see \cite{Nielsen} Theorem 3.2, generalizing \cite{Segal} Prop. 4.1 in topological case).
Now let $t_0\in T$ be an element of finite order $N$.
Specializing \eqref{loc-map1} with respect to the homomorphism $\Tr(t_0,?):R\to\Z[\sqrt[N]{1}]$
we get a map
\begin{equation}\label{loc-map}
K_{t_0}=K^T_0(X)\otimes_R\Z[\sqrt[N]{1}]\to\bigoplus_{p\in X^T}\Z[\sqrt[N]{1}]
\end{equation}

\begin{lem}\label{inj-lem}
Assume that $K^T_0(X)$ is a free $R$-module, and
$\det(1-t_0,T_pX)\neq 0$ for every $p\in X^T$. 
Then the map \eqref{loc-map} is injective and 
becomes an isomorphism after tensoring with $\Q$. 
\end{lem}

\Pf . First, Theorem 3.2 of \cite{Nielsen} easily implies that in the case when $K^T_0(X)$ is a free
$R$-module the map \eqref{loc-map1} becomes an isomorphism after
inverting in $R$ all elements $\chi-1$, where 
$\chi$ is a character 
occuring in the $T$-action on one of the tangent spaces $T_pX$ for $p\in X^T$.
Therefore, under our assumptions on $t_0$ the map \eqref{loc-map} becomes an isomorphism
after tensoring with $\Q$. Since it is a map of free $\Z[\sqrt[N]{1}]$-modules, the injectivity follows.
\ed

For a class $F\in K^T_0(X)$ we denote by $v(F)$ the corresponding class
in $K_{t_0}$. We also set
$v(F)_p=\Tr(t_0,F|_p)\in \Z[\sqrt[N]{1}]$ for $p\in X^T$, so that
the map \eqref{loc-map} sends $v(F)$ to $(v(F)_p)_{p\in X^T}$.

\begin{lem}\label{line-lem} 
Under the assumptions of Lemma \ref{inj-lem} for a $T$-equivariant line bundle $L$ on $X$
the class $v(L^m)\in K_{t_0}$ depends only on the remainder $m\mod(N)$.
\end{lem}

\Pf . For each $p\in X^T$ let $\chi_p$ be the character of the torus $T$ corresponding to its action
on $L|_p$. Then $v(L^m)_p=\chi_p(t_0)^m$ depends only on the remainder $m\mod(N)$.
Now the assertion follows from the injectivity of \eqref{loc-map} in this case.
\ed

For a $T$-equivariant vector bundle $V$ on $X$ (and hence, for an object of $D(X/T)$)
we have the following Lefschetz type formula for the equivariant
Euler characteristic:
\begin{equation}\label{Lefschetz-eq1}
\chi^T(X,V)=\sum_{p\in X^T}\det(1-t,T_pX)^{-1}[V|_p],
\end{equation}
where $[V|_p]\in R$ (see \cite{Nielsen}, 4.9).
Hence, for $t_0$ such that $\det(1-t_0,T_p(X))\neq 0$ for all $p\in X^T$,
for classes $[V],[W]\in K^T_0(X)$ one has
\begin{equation}\label{Lefschetz-eq2}
\Tr(t_0,\chi^T(V,W])=\sum_{p\in X^T}c_p v(V^*)_p v(W)_p,
\end{equation}
where $c_p=\det(1-t_0,T_pX)^{-1}$.
Since the eigenvalues of $t_0$ are roots of unity, we have 
$v(V^*)_p=\ov{v(V)_p}$ (where bar denotes complex conjugation).
Let us equip the free $\Z[\sqrt[N]{1}]$-module $\oplus_{p\in X^T} \Z[\sqrt[N]{1}]$
with the sesquilinear form 
$$h(v,w)=\sum_{p\in X^T}c_p \ov{v}_pw_p,$$
and let $H$ denote the pull-back of this form to $K_{t_0}$ via \eqref{loc-map}.
Then we can rewrite \eqref{Lefschetz-eq2} as
\begin{equation}\label{Tr-H-eq}
\Tr(t_0,\chi^T(V,W))=H(v(V),v(W)).
\end{equation}

Recall that by Serre duality we have
$$\chi^T(V,W)^*=(-1)^{\dim X}\chi^T(W,V\otimes\om_X).$$
Taking traces of $t_0$ and using \eqref{Tr-H-eq} we obtain
$$\ov{H(v(V),v(W))}=(-1)^{\dim X}H(v(W),v(V\otimes\om_X)).$$
The condition $(\star)$ implies that $t_0$ acts as $(-1)^{\dim X}$ on the fibers
of $\om_X$ at all points $p\in X^T$.
Hence, $v(V\otimes\om_X)_p=(-1)^{\dim X}v(V)_p$ for all $p\in X^T$, and the above equation becomes
$$\ov{H(v(V),v(W))}=H(v(W),v(V)),$$
which implies part (i) of Theorem \ref{main-thm} because of \eqref{Tr-H-eq}.

Now part (ii) follows from the observation that if $E$ is exceptional then
$\chi^T(E,E)=1$ (resp, if $(E_1,E_2)$ is an exceptional pair then $\chi^T(E_2,E_1)=0$).

Part (iii) is an immediate consequence of the following number theoretic result (with $r=1$).

\begin{prop}\label{number-prop} 
Let $z_1,\ldots,z_n\in\Z[\sqrt[N]{1}]$ be such that
$|z_1|^2+\ldots+|z_n|^2=r$, where $r$ is a rational number, $0<r\le 1$. 
Then for some $i_0$ we have $z_i=0$ for $i\neq i_0$, and $z_{i_0}$ is a root of unity (so $r=1$).
\end{prop}

\Pf . Since the Galois group of $\Q(\sqrt[N]{1})$ over $\Q$ is abelian, for every element
$\si$ in this group we have $\si(|z_i|^2)=|\si(z_i)|^2$. Hence, applying $\si$ to the equality
$|z_1|^2+\ldots+|z_n|^2=r$ we get that all conjugates of $|z_i|^2$ are positive real numbers
$\le 1$ for each $i$. Assume that $z_1\neq 0$. Then the fact that the norm of $|z_1|^2$ is an integer
implies that $|z_1|^2=1$. Thus, all conjugates of $z_1$ have absolute value $1$. 
By Kronecker's theorem, it follows that $z_1$ is a root of unity.
\ed

\begin{rem} The above proof works also in a more general situation when $z_i$'s are integers
in some Galois extension of $\Q$ such that the complex conjugation induces a central element
of the corresponding Galois group.
\end{rem}

To prove part (iv) of Theorem \ref{main-thm} we start by observing 
that the action of $N$ on $X$ preserves $X^T$. Furthermore, for
$p\in X^T$ and an $N$-equivariant object $E$ we have
\begin{equation}\label{conj-wp-eq}
\ov{v(E)_{wp}}=\Tr(t_0^{-1},E|_{wp})=\Tr(wt_0w^{-1}, E|_{wp})=\Tr(t_0, E|_p)=v(E)_p.
\end{equation}
From (iii) we derive that if $E$ is exceptional then $v(E)=\pm\zeta v(E_i)$ for some $i$.
This implies that $v(E)_p=\pm\zeta v(E_i)_p$ for all $p\in X^T$.
It remains to observe that \eqref{conj-wp-eq} applies both to $E$ and $E_i$ and that the map
\eqref{loc-map} is injective by Lemma \ref{inj-lem}
(recall that by Lemma \ref{exc-basis-lem}, $K_0^T(X)$ is
a free $R$-module). Hence, the coefficient of proportionality between $v(E)$ and $v(E_i)$
should be real, so $v(E)=\pm v(E_i)$.

\medskip

\noindent
{\it Proof of Corollary \ref{rk-cor}}. 
Let $\Phi_N$ denote the $N$th cyclotomic polynomial. Since for $N=p^k$ one has
$\Phi_N(1)\equiv 0\mod(p)$, there is a well-defined ring homomorphism
$\rho:\Z[\sqrt[N]{1}]\to\Z/p\Z$ sending $\zeta_N$ to $1$. The composition
$\rho\circ\Tr(t_0,?):R\to\Z/p\Z$ coincides with the reduction modulo $p$ of the homomorphism
$\Tr(1,?):R\to\Z$. 
Therefore, from Lemma \ref{eval-lem} we deduce an isomorphism
$$K_{t_0}\otimes_{\Z[\sqrt[N]{1}]}\Z/p\Z\simeq K_0(X)\otimes\Z/p\Z$$
compatible with the Euler forms. 
It follows that $\chi\mod(p)$ is symmetric, and for any object $F$ of $D(X)$
equipped with $T$-equivariant structure, the class associated with $F$ in $K_0(X)\otimes\Z/p\Z$
is obtained from the class $v(F)\in K_{t_0}$ using the specialization
with respect to $\rho$. This immediately implies the first assertion. 
Note that every exceptional object $E\in D(X)$ admits a $T$-equivariant
structure (see Lemma \ref{equiv-lem}). Hence,
by part (iii) of Theorem \ref{main-thm}, we have $v(E)=\pm\zeta_N^i\cdot v_i$ for some $i$.
Applying the homomorphism $\rho$ we deduce \eqref{mod-p-congr}, which in turn implies
\eqref{sum-chi-congr}.
\ed

\section{Central equivariant objects and applications}\label{appl-sec}

The following result due to Merkurjev is explained as Lemma 2.9 in \cite{VV} (it is a combination
of Cor. 2.15 and Prop. 4.1 of \cite{Mer}).

\begin{lem}\label{Mer-lem}
Assume that $T$ is a maximal torus in a connected reductive group $G$ 
such that the commutant of $G$ is simply connected, and let $X$ be a smooth projective variety with an action of $G$. Then the natural morphism
$$K^G(X)\otimes_{R(G)}R(T)\to K^T(X)$$
is an isomorphism.
\end{lem}

\begin{defi} Let $X$ be a variety equipped with an action of an algebraic group $G$. Assume that the center $Z_G\sub G$ acts trivially on $X$. 
Then we have a decomposition of the category $\Coh^G(X)$ into
the direct sum of subcategories $\Coh^G(X)_{\chi}$, where $\chi$ runs through characters of $Z_G$
and $Z_G$ acts via $\chi$ on $G$-equivariant coherent sheaves in $\Coh^G(X)_{\chi}$.
We say that an object $V\in D(X/G)$ is {\it central} if it belongs to $D(\Coh^G(X)_{\chi})$ for some
$\chi:Z_G\to\G_m$ (called the {\it central character} of $V$).
\end{defi}

For example, any indecomposable object in $D(X/G)$ is central. 
Hence, every exceptional object equipped with a $G$-equivariant structure is central.
The tensor product of central objects is again central (and the central characters get multiplied).

Using central $G$-equivariant bundles we get a simple way of decomposing
$K_{t_0}=K^T_0(X)\otimes_R \Z[\sqrt[N]{1}]$ into $\Z[\sqrt[N]{1}]$-submodules,
mutually orthogonal with respect to the specialization of the
equivariant Euler form $\Tr(t_0,\chi^T(\cdot,\cdot))$.

\begin{prop}\label{central-orth-prop} Let $G$ be a connected
reductive group with simply connected commutant,
and let $X$ be a smooth projective variety with an action of $G/Z_G$. Let $t_0\in T$ be an element
of order $N$ in a maximal torus in $G$. 
Assume that $K^T_0(X)$ is a free $R$-module, and
$\det(1-t_0,T_pX)\neq 0$ for every $p\in X^T$. 
Assume also that for some element $z_0\in Z_G$ and some element
$w_1\in W$ in the Weyl group of $G$ one has $w_1(t_0)=z_0t_0$. 
Then we have an orthogonal (with respect to the Euler form $\Tr(t_0,\chi^T(\cdot,\cdot))$) direct sum decomposition:
$$K_{t_0}=\bigoplus_{\zeta^n=1}K_{t_0}(\zeta),$$
where $n$ is the order of $z_0$, and for each $n$-th root of unity $\zeta$ we denote by $K_{t_0}(\zeta)$ the subgroup
in $K_{t_0}$ generated by the classes of $G$-equivariant bundles on which $z_0$ acts by $\zeta$.
If in addition, $N$ is a power of a prime $p$ then
we have a similar decomposition of $K_0(X)\otimes\Z/p\Z$, orthogonal with respect to the
reduction of the Euler form.
\end{prop}

\Pf .
By Lemma \ref{Mer-lem}, the classes of $G$-equivariant bundles generate
$K^T_0(X)$ over $R$. Also, Lemma \ref{inj-lem} together with \eqref{Tr-H-eq} imply that
the bilinear form $\Tr(t_0,\chi^T(\cdot,\cdot))$ on $K_{t_0}$ is nondegenerate.
Hence, it suffices to check that the pieces $K_{t_0}(\zeta)$ and
$K_{t_0}(\zeta')$ are orthogonal for $\zeta\neq\zeta'$. 
It is enough to check that if $V$ is a $G$-equivariant vector bundle
such that $z_0$ acts on $V$ by a scalar $\zeta\neq 1$ 
then $\Tr(t_0,\chi^T(V))=0$. To this end we observe
that $\chi^T(V)$ is a $W$-invariant element in $R$. Furthermore, $\chi^T(V)$ belongs
to the $\Z$-span of characters $\chi$ of $T$ such that $\chi(z_0)=\zeta$.
Therefore, we get
$$\Tr(t_0,\chi^T(V))=\Tr(w_1(t_0),\chi^T(V))=\Tr(z_0t_0,\chi^T(V))=\zeta\Tr(t_0,\chi^T(V))$$
which implies the required vanishing.
\ed


We will combine Proposition \ref{central-orth-prop} with the following result employing the same idea
as in part (iv) of Theorem \ref{main-thm}.

\begin{prop}\label{Galois-prop}
Let $G$ be a connected reductive group with simply connected commutant, 
and let $X$ be a smooth projective variety with an action of $G/Z_G$. Let $t_0\in T$ be an element
of order $N$ in a maximal torus in $G$. We assume that $T$ comes from a split torus over $\Q$ and
we use the corresponding Galois action on the elements of finite order in $T$.
Assume that for every $\si\in\Gal(\Q(\sqrt[N]{1})/\Q)$
there exists $w_{\si}\in W$ such that $\si(t_0)=w_{\si}(t_0)$.
Let us denote by $M\sub K_{t_0}$
the $\Z$-span of the classes of $G$-equivariant bundles on $X$.
Then the natural ring homomorphism
$$M\otimes_{\Z}\Z[\sqrt[N]{1}]\to K_{t_0}$$
is an isomorphism. If $(E_1,\ldots,E_k)$ is a full exceptional collection in $D(X)$, where each
$E_i$ can be equipped with $G$-equivariant structure, then $M$ coincides with the $\Z$-span
of the classes of $E_i$'s.
\end{prop}

\Pf . By Lemma \ref{Mer-lem}, we have
$$K_{t_0}\simeq K_0^G(X)\otimes_{R(G)}\Z[\sqrt[N]{1}],$$
where we use the ring homomorphism $R(G)=R(T)^W\to\Z[\sqrt[N]{1}]$ 
induced by $\Tr(t_0,?):R(T)\to\Z[\sqrt[N]{1}]$. It remains to note that our assumptions imply that for every
$W$-invariant element $f\in R(T)$ the element $\Tr(t_0,f)\in \Z[\sqrt[N]{1}]$ 
will be invariant with respect to
$\Gal(\Q(\sqrt[N]{1})/\Q)$, hence, an integer. For the last assertion one has to use the fact that
the classes of $E_i$'s form a basis of $K_0^G(X)$ over $R(G)$ (see Lemma \ref{exc-basis-lem}).
\ed

Now let us specialize to the particular case $X=G/P$, where $G$ is a simply connected semisimple group, $P$ is a maximal parabolic subgroup of $G$.
Let $\De$ denote the set of roots associated with $G$, and let
$\Pi=(\a_1,\ldots,\a_n)$ denote the set of simple roots. 
Assume that $P$ is the standard maximal parabolic in $G$ associated with a simple root $\a_i$ 
(or rather, with the subset $\Pi\setminus\{\a_i\}\sub\Pi$). 
Then the weights of the maximal torus $T$ on 
${\frak g}/{\frak p}$ are exactly $-\a$, where $\a$ is a positive root
in which $\a_i$ enters with nonzero multiplicity, i.e., $(\a,\om_i)>0$, where
 $\om_i$ is the fundamental weight associated with $\a_i$.
The set of $T$-stable points
$X^T$ is in bijection with $W/W_P$, where $W=N/T$ is the Weyl group of $G$, and
$W_P\sub W$ is the subgroup
generated by the reflections with respect to $\Pi\setminus\{\a_i\}$. 
If $p(w)\in X^T$ denotes the point corresponding to $wW_P\in W/W_P$ then 
the weights of $T$ on the tangent space $T_{p(w)}X$ are $-w\a$, 
where $\a$ is as above. 
Thus, the assumption $(\star)$ in Theorem \ref{main-thm} can be reformulated
in the following form.

\begin{prop}\label{Weyl-lem2}
Let $G$ be a simply connected semisimple group, 
and let $X=G/P$, where $P$ is the standard maximal
parabolic subgroup associated with a simple root $\a_i$.
Set 
$$N=N_i=\sum_{\a\in\De:(\a,\om_i)>0}\frac{(\a,\om_i)}{(\om_i,\om_i)}.$$ 
Then an element $t_0$ of the maximal torus $T$ satisfies the assumption $(\star)$ of Theorem
\ref{main-thm} if and only if the following two conditions are satisfied:

\noindent
(a) for every root $\a$ one has $\a(t_0)\neq 1$; 

\noindent
(b) $\om_i(t_0)^N=(-1)^{\dim X}$ and for every
$\a\in\De$ such that $||\a||=||\a_i||$ one has $\a(t_0)^N=1$.
\end{prop}

\Pf . The set of weights of $T$ on ${\frak g}/{\frak p}$ is invariant with respect to $W_P\sub W$. 
Hence, the sum of these
weights is proportional to $\om_i$. It follows that this sum equals $-N\om_i$, so we have
$$\det(t,{\frak g}/{\frak p})=\om_i(t)^{-N}.$$
Thus, the condition $(\star)$ can be restated as follows:
\begin{equation}\label{w-om-eq}
(w\om_i)(t_0)^N=(-1)^{\dim X}
\end{equation} 
and $(w\a)(t_0)\neq 1$ for all $w\in W$ and all roots $\a$ such that $(\a,\om_i)\neq 0$. 
It is easy to see that in fact every root can be written in the form
$w\a$ with $(\a,\om_i)\neq 0$ for some $w\in W$. Hence, the inequalities are equivalent to the condition
$\a(t_0)\neq 1$ for all roots $\a$. On the other hand, since $s_i\om_i=\om_i-\a_i$,
the equalities \eqref{w-om-eq} imply that $\la(t_0)^N=1$ for all $\la$ in the lattice $Q_i$
spanned by $(w\a_i)_{w\in W}$ (i.e., by the set of roots of the same length as $\a_i$). 
Note that since $s_j\om_i=\om_i-\de_{ij}\a_j$,
the $Q_i$-coset $\om_i+Q_i$ is stable under $W$. 
This implies our statement.
\ed

\begin{cor}\label{odd-cor}
In the situation of the above Proposition, assume that $\dim X$ is odd. Then 
the necessary condition for the existence of an element
$t_0$ satisfying $(\star)$ is that the class of $\om_i\mod(Q_i)$ has an even order,
where $Q_i$ is the sublattice of the weight lattice spanned by all roots of the same length as $\a_i$. 
\end{cor}

Using Proposition \ref{Weyl-lem2} we will be able to check for almost all of the spaces $G/P$ (where $P$ is a maximal parabolic) whether an element $t_0$ satisfying $(\star)$ exists, leaving out several
cases of types $E_7$ and $E_8$. 
In the next Proposition we use notations from the Tables I-IX in \cite{Bou}. 
Recall that since $G$ is simply connected, the
character lattice of the maximal torus $T$
coincides with the weight lattice of the corresponding root system.
For example, for type $C_n$ this gives an identification of $T$ with $\G_m^n$.
For types $B_n$ and $D_n$ the weight lattice is spanned by the
standard lattice $\Z^n$ together with the vector $(\sum_{i=1}^n \vareps_i)/2$. Hence, we obtain
$$T=\{(x_1,\ldots,x_n;x)\in\G_m^{n+1}\ |\ x^2=\prod x_i\}$$
in these cases. The maximal torus for type $F_4$ has the same description with $n=4$.

\begin{prop}\label{cases-prop}
Let $X=G/P$, where $G$ is a simply connected simple algebraic group, and $P$ is the standard
maximal parabolic subgroup associated with the simple root $\a_i$. 

\noindent
(a) If $G$ is of classical type then $t_0$ satisfying $(\star)$ exists only in the following cases:

(i) $G$ is of type $A_n$, arbitrary $i$; 

(ii) $G$ is of type $B_n$, $i=1$ or $i=n$;

(iii) $G$ is of type $C_n$, $i=1$ or $i=n$;

(iv) $G$ is of type $D_n$, $i=1,n-1$, or $n$.

In the cases (ii)--(iv) $t_0$ can be chosen in such a way that 
\begin{equation}\label{w0-t0-eq}
w_0(t_0)=t_0^{-1},
\end{equation} 
where $w_0\in W$ is the element of maximal
length.

\noindent
(b) $t_0$ satisfying $(\star)$ does not exist if $G$ is of type $G_2$ or $F_4$.

\noindent
(c) If $G$ is of type $E_6$ then $t_0$ satisfying $(\star)$ exists if and only if either $i=1$ or $i=6$.
In these cases it can be chosen in such a way that \eqref{w0-t0-eq} holds.

\noindent
(d) If $G$ is of type $E_7$ and $i=1,3$, or $4$ (resp., $G$ is of type $E_8$ and $i=6, 7$ or $8$)
then $t_0$ satisfying $(\star)$ does not exist.

\noindent
(e) if $G$ is of type $E_7$ and $i=7$
then $t_0$ satisfying $(\star)$ exists and \eqref{w0-t0-eq} holds.
\end{prop}

\Pf . (a) For type $A_n$ the maximal torus is 
$T=\{(x_1,\ldots,x_{n+1})\in\G_m^{n+1}\ |\ \prod x_i=1\}$. Also,
with the notations of Proposition \ref{Weyl-lem2} we have $N_i=n+1$. Thus, the conditions
of this Proposition for the element $t_0=(x_1,\ldots,x_{n+1})$ become
$$(x_1\ldots x_i)^{(n+1)}=(-1)^{i(n+1-i)},$$ 
$$x_k^{n+1}=x_l^{n+1} \text{ and } x_k\neq x_l \text{ for }k<l.$$
Thus, these conditions are satisfied when $\{x_1,\ldots,x_{n+1}\}$ is the set of all $(n+1)$th
roots of $(-1)^{i(n+1-i)}$. 

For type $B_n$ the element $t_0=(x_1,\ldots,x_n;x)$ should satisfy $x_k\neq 1$ for all $k$
and $x_k\neq x_l^{\pm 1}$ for $k<l$. For $i<n$ we have $\om_i=\vareps_1+\ldots+\vareps_i$,
$N_i=2n-i$, $\dim X=i(i+1)/2+2i(n-i)$, and
the lattice $Q_i\sub \Z^n$ consists of all vectors with even sums of coordinates.
Thus, we should have $x_k^{2N_i}=1$ for all $k$, $x_k^{N_i}=x_l^{N_i}$ for $k<l$.
In other words, $x_k^{N_i}=\pm 1$ does not depend on $k$, which is impossible for $i>1$.
In the case $i=1$ all $x_k$'s should be $(2n-1)$th roots of $-1$, and we can set $x_1=-1$ and
choose $\{x_2,\ldots,x_{n}\}$ to contain one from each conjugate pair of the remaining $2n-2$ roots
(and let $x$ be a square root of $\prod x_k$).
In the case $i=n$ we have $\om_n=(\sum \vareps_k)/2$, 
$N_n=2n$, $\dim X=n(n+1)/2$, and $Q_i=\Z^n$. Thus, $x_k$'s should
be $(2n)$th roots of $1$, and $x$ should satisfy $x^{2n}=(-1)^{n(n+1)/2}$. Hence, we can set 
$x_k=\zeta_{2n}^k$ for $k=1,\ldots,n$, and let $x$ be a square root of $\prod x_k$.
The condition \eqref{w0-t0-eq} for types $B_n$ and $C_n$ is automatic since $w_0$ sends
every $t\in T$ to $t^{-1}$.

For type $C_n$ the element $t_0=(x_1,\ldots,x_n)$ should satisfy $x_k\neq\pm 1$ for all $k$
and $x_k\neq x_l^{\pm 1}$ for $k<l$. We have $\om_i=\vareps_1+\ldots+\vareps_i$, $N_i=2n-i+1$,
$\dim X=i(i+1)/2+2i(n-i)$. For $i<n$ the lattice $Q_i\sub\Z^n$ consists of all vectors with even sums
of coordinates. Hence, as in the case of $B_n$ we deduce that $x_k^{N_i}=\pm 1$ does not depend
on $k$ which is impossible for $i>1$ (recall that now we have the condition $x_k\neq\pm1$).
In the case $i=1$ we can choose $\{x_1,\ldots,x_n\}$ to contain one from each conjugate pair
of $(2n)$th roots of $-1$. In the case $i=n$ the lattice $Q_n\sub\Z^n$ consists of all vectors with
even coordinates. Hence, we can set $x_k=\zeta_{2n+2}^k$ for $k=1,\ldots,n$.

For type $D_n$ the element $t_0=(x_1,\ldots,x_n;x)$ should satisfy $x_k\neq x_l^{\pm 1}$ for $k<l$.
For $i<n-1$ we have $\om_i=\vareps_1+\ldots+\vareps_i$, $N_i=2n-i-1$, $\dim X=i(i-1)/2+2i(n-i)$, the lattice
$Q_i=Q\sub\Z^n$ consists of all vectors with even sums of coordinates. As above we can rule
out the cases $1<i<n-1$. For $i=1$ we can take $x_k=\zeta_{2n-2}^{n-k}$, $k=1,\ldots,n$
(and let $x$ be a square root of $\prod x_k$). For $i=n$ we have $\om_n=(\sum \vareps_k)/2$,
$N_n=2n-2$, $\dim X=n(n-1)/2$. Hence, we can use the same element $t_0$ as for $i=1$.
If $n$ is even then the condition \eqref{w0-t0-eq} is automatic. If $n$ is odd then
$$w_0(x_1,\ldots,x_n;x)=(x_1^{-1},\ldots,x_{n-1}^{-1},x_n;x^{-1}x_n),$$
so the condition \eqref{w0-t0-eq} holds since $x_n=1$.
The case $i=n-1$ follows by the symmetry of the Dynkin diagram.

\noindent
(b) The case of type $G_2$ and the cases $i=1,4$ of type $F_4$ follow 
immediately from Corollary \ref{odd-cor}. In the remaining two cases $i=2,3$ for type $F_4$ we have $N_2=5$, $N_3=7$, and the element $t_0=(x_1,\ldots,x_4;x)$ should satisfy $x_k\neq 1$, 
$x_k\neq x_l^{\pm 1}$ for $k<l$, and
$x_k^{N_i}=\pm 1$ does not depend on $k$, which is impossible.

\noindent
(cd) The non-existence of $t_0$ in all the relevant cases follows Corollary \ref{odd-cor}.
It remains to consider the case of type $E_6$ and $i=1$
(the case $i=6$ will follow by the symmetry of the Dynkin diagram). We have
$N_1=12$. Thus, the element $t_0$ should satisfy $t_0^{12}=1$ and $\a(t_0)\neq 1$
for every root $\a$. Let $\La$ denote the weight lattice. Then the group of elements
of order $12$ in $T$ is canonically dual to the finite group $\La\otimes\Z/12\Z$.
Thus, to give $t_0$ it is enough to specify an element of order $12$ in $\La^{\vee}\otimes\Q/\Z$.
Equivalently, we have to produce an element $\la\in\La\otimes\Q$ such that
$12(\la,\om)\in\Z$ for every weight $\om$ and $(\la,\a)\not\in\Z$ for every root $\a$.
Set $12\la=\sum_{i=1}^5 a_i \vareps_i+b(\vareps_8-\vareps_7-\vareps_6)$. Then the conditions can be rewritten in
terms of these coordinates as follows: $a_i\in \frac{1}{2}\Z$, $a_i\equiv a_j\mod(\Z)$ for all $i,j$,
$c:=\frac{1}{2}(3b+\sum_{i=1}^5 a_i)\in\Z$ and 
$$a_i\not\equiv \pm a_j\mod(12\Z) \text{ for } i\neq j,$$
$$c\not\equiv\sum_{i\in S}a_i\mod(12\Z) \text{ for } S\sub[1,5] \text{ with } |S| \text{ even}.$$
If we want in addition to have $w_0(t_0)=t_0^{-1}$ then we should impose two relations:
$a_2-a_1=a_4-a_3$ and $c=a_2+a_3+2a_5$. It is easy to check that
$$(a_1,\ldots,a_5;c)=(0,1,2,3,4;11)$$
is a solution.

\noindent
(e) The condition \eqref{w0-t0-eq} holds automatically in this case since $w_0$ sends every
$t\in T$ to $t^{-1}$. We have $N_7=18$ and $\dim X=27$, so to give
$t_0$ satisfying $(\star)$ is equivalent to finding a rational weight $\la$ such that
$18(\la,\a)\in\Z$, $(\la,\a)\not\in\Z$ for every root $\a$, while $18(\la,\om_7)-1/2\in\Z$.
Set $18\la=\sum_{i=1}^6 a_i \vareps_i+b(\vareps_8-\vareps_7)$. Then we should have
$a_i\in \frac{1}{2}\Z$, $a_i\equiv a_j\mod(\Z)$ for all $i,j$, $\sum_{i=1}^6 a_i\equiv 1\mod(2\Z)$, 
$c:=b+\frac{1}{2}\sum_{i=1}^6 a_i\equiv a_6\mod(\Z)$, and
$$a_i\not\equiv \pm a_j\mod(18\Z) \text{ for } i\neq j, \ \ 2c\not\equiv\sum_{i=1}^6 a_i\mod(18\Z),$$
$$c\not\equiv\sum_{i\in S}a_i\mod(18\Z) \text{ for } S\sub[1,6] \text{ with } |S| \text{ odd}.$$ 
We can take
$$(a_1,\ldots,a_6;c)=(0,1,2,3,4,5;16) \text{ or }(0,1,2,3,4,5;17)$$
as solutions. 
\ed

Now let us consider applications of Theorem \ref{main-thm} (and of Propositions
\ref{central-orth-prop} and \ref{Galois-prop}) to concrete varieties.

\bigskip

{\it Projective spaces}

\medskip

In this case it is more convenient to work with the action of $\GL_n$ on $\P^{n-1}$, rather than
$\SL_{n}$, so $T$ will denote the set of diagonal matrices in $\GL_n$.
Let $p_1,\ldots,p_{n}\sub\P^{n-1}$ denote the $T$-fixed points where the action of $T$ on
the fiber of $\OO_{\P^{n-1}}(1)$ at $p_i$ corresponds to the $i$-th coordinate character
$\vareps_i:T\to\G_m$.

As in Proposition \ref{cases-prop} one checks that the element
\begin{equation}\label{GL-t0-eq}
t_0=(1,\zeta_{n},\zeta_{n}^2,\ldots,\zeta_{n}^{n-1})
\end{equation}
satisfies the assumption $(\star)$ of Theorem \ref{main-thm} (where $\zeta_{n}$ is
a primitive $n$-th root of unity).

Let $V$ be a central object in $D(\P^{n-1}/\GL_{n})$.
The action of the center $Z_{\GL_{n}}=\G_m$ on $V$ 
is given by the character $\G_m\to\G_m:\la\mapsto\la^{m}$ for some $m\in\Z$. 
Let us denote by $m(V)\in \Z/n\Z$ the remainder of $m$ modulo $n$.
Note that if we tensor $V$ with a character of $\GL_{n}$ then $m(V)$ will not change.
In particular, for an exceptional object $V$ in $D(\P^{n-1})$ we can define $m(V)\in\Z/n\Z$
by choosing any $\GL_n$-equivariant structure on $V$ (see Lemma \ref{equiv-lem}).
Note that the center of $\GL_n$ acts on $\OO_{\P^{n-1}}(1)$ through the identity character.
This implies that for a central object $V$ on $D(\P^{n-1}/\GL_{n})$ one has
$$m(V)\rk(V)\equiv \deg(V)\mod(n).$$

Let $w_1\in W=S_{n}\sub\GL_{n}$ be the cyclic permutation such that $w_1(p_i)=p_{i-1}$.
Then we have 
$$w_1(t_0)= \zeta_{n}\cdot t_0,$$
where $\zeta_{n}$ is viewed a scalar matrix in $T$.
Thus, the assumptions of Propositions \ref{central-orth-prop} and
Proposition \ref{Galois-prop} are satisfied in this case: for the latter
we can take $w_{\si}\in S_{n}$ to be
the permutation of the set of $n$th roots of unity induced by 
$\si\in \Gal(\Q(\sqrt[n]{1})/\Q)$. 
Hence, we derive the following corollary from these Propositions. 

\begin{cor}\label{equiv-proj-cor} 
Let $V$ be a central object in $D(\P^{n-1}/\GL_{n})$. Then
for some integer $a(V)$ one has 
\begin{equation}\label{aV-eq}
v(V)=a(V)v(\OO(m(V)))
\end{equation}
in $K^T_0(\P^{n-1})\otimes_R\Z[\sqrt[n]{1}]$ (recall that $v(\OO(m))$ depends only on
$m\mod(n)$ by Lemma \ref{line-lem}).
\end{cor}

\noindent
{\it Proof of Theorem \ref{congr-proj-thm}.}
In the case when $n=p^k$, where $p$ is prime, we can apply the homomorphism
$\rho:\Z[\sqrt[n]{1}]\to\Z/p\Z$ that sends $p^k$th roots of unity to $1$. Then
we get from Corollary \ref{equiv-proj-cor} the following congruence in $K_0(\P^{n-1})\otimes\Z/p\Z$:
$$[V]\equiv a(V)[\OO(m(V))].$$
Taking ranks of both sides we see that $a(V)\equiv \rk(V)\mod(p)$.
This immediately implies part (i) of Theorem \ref{congr-proj-thm}.
For part (ii) we have to recall that by Lemma \ref{equiv-lem} an exceptional object $E$ admits
a $\GL_n$-equivariant structure, unique up to tensoring with a character of $\GL_n$. Also,
by Theorem \ref{main-thm}(iii) we see that $a(E)=\pm 1$.
Finally, part (iii) follows immediately from Corollary \ref{rk-cor}.
\ed

\begin{rem} The fact that $a(E)=\pm 1$ for an exceptional object $E$ can also be proven
by calculating $1=\Tr(t_0,\chi^T(E,E))$ using \eqref{aV-eq}. Also, part (iii) of Theorem
\ref{congr-proj-thm} can be deduced from a version of \eqref{aV-eq} that holds for products
of projective spaces. This would make a proof of Theorem \ref{congr-proj-thm} completely
independent from Theorem \ref{main-thm}. Because of this it is possible to generalize
this argument to relative projective spaces, see Theorem \ref{rel-proj-thm} below.
\end{rem}

\bigskip

{\it Grassmannians}

\medskip

Let $T$ be the maximal torus of $\GL_n$ acting on the Grassmannian $G(k,n)$ of $k$-planes in
the $n$-dimensional space in the standard
way. It is easy to see that the element \eqref{GL-t0-eq} still satisfies the condition $(\star)$
of Theorem \ref{main-thm} (cf. Proposition \ref{cases-prop}).

As before, for a central object
$V$ in $D(G(k,n)/\GL_n)$ with the central character $\la\mapsto\la^m$ 
we set $m(V)=m\mod(n)$. This remainder does not change under tensoring
$V$ with a character of $\GL_n$.
We denote by $\OO(1)$ the ample generator of $\Pic(G(k,n))$. Note that
it has the central character $\la\mapsto\la^k$. Therefore, one has
\begin{equation}\label{m-deg-eq}
m(V)\rk(V)\equiv k\deg(V)\mod(n).
\end{equation}
 


Recall that if $\UU$ is the tautological rank $k$ bundle on $G(k,n)$ then
the vector bundles $\Sigma^{\la}\UU$, where $\la$ runs over partitions 
$n-k\ge\la_1\ge\ldots\ge\la_k\ge 0$,
can be ordered to form a full exceptional collection on $G(k,n)$ (see \cite{Kap1}).
Hence, the classes $v(\Sigma^{\la}\UU)$ with $\la$ as above form a basis of
$K_0^T(G(k,n))\otimes_R\Z[\sqrt[n]{1}]$ over $\Z[\sqrt[n]{1}]$.

As in the case of projective spaces, using Propositions \ref{central-orth-prop}
and \ref{Galois-prop} we derive the following.

\begin{cor}\label{equiv-Grass-cor} 
(i) For a central object $V$ in $D(G(k,n)/\GL_n)$ the class 
$v(V)\in K_0^T(G(k,n))\otimes_R\Z[\sqrt[n]{1}]$
is a linear combination with integer coefficients of the classes $v(\Sigma^{\la}\UU)$, where
$\la$ runs over partitions $n-k\ge\la_1\ge\ldots\ge\la_k\ge 0$ such that
$\la_1+\ldots+\la_k\equiv m(V)\mod(n)$.

\noindent
(ii) Let $V$ and $V'$ be central objects in $D(G(k,n)/\GL_n)$ with
$m(V)\not\equiv m(V')\mod(n)$.
Assume that $n=p^r$, where $p$ is a prime. Then $\chi(V,V')\equiv 0\mod(p)$.
\end{cor}

\begin{rem}
The (integer) structure constants of the multiplication on $K_0^T(G(k,n))\otimes_R\Z[\sqrt[n]{1}]$
can be easily computed
from the Littlewood-Richardson rule. One just has to observe (looking at the definition of the Schur
functions) that for an arbitrary partition $\la_1\ge\ldots\ge\la_k\ge 0$ the class $v(\Sigma^{\la}\UU)$ up to a sign depends only on the residues modulo $n$ of the numbers $\la_1+k-1,\la_2+k-2,\ldots,\la_k$. More precisely, if for some $i\neq j$ we have $\la_i-i\equiv\la_j-j\mod(n)$ then $v(\Sigma^{\la}\UU)=0$. Otherwise, the class $v(\Sigma^{\la}\UU)$ coincides up to a sign with one of the basis classes.
\end{rem}

\begin{thm}\label{Grass-thm} 
(i) Let $(E_1,\ldots, E_s)$ be a full exceptional collection 
in $D(G(k,n))$ (so that $s={n\choose k}$).
Assume that $k$ is relatively prime to $n$. Then for each $m\in\Z/n\Z$ exactly
$s/n$ objects $E_i$ from the collection have $m(E_i)\equiv m\mod(n)$.

\vspace{2mm}

\noindent
(ii) Let $p$ be a prime, and let $1\le k\le p-1$.
For every exceptional object $E$ in $D(G(k,p))$ one has $\rk(E)\not\equiv 0 \mod(p)$.
The same congruence holds for exceptional objects on the products
$G(k_1,p)\times\ldots\times G(k_l,p)$. If $(E_1,\ldots,E_s)$ is a full exceptional collection in $D(G(k,p))$ 
then for every $\mu\in\Z/p\Z$ exactly $s/p$ of the objects $E_i$ have
$\deg(E_i)/\rk(E_i)\equiv\mu\mod(p)$.

\vspace{2mm}

\noindent
(iii) Now let $n=p^r$, where $p$ is prime, and let $V$ be a central object in $D(G(p,n)/\GL_n)$ with 
$m(V)\not\equiv 0\mod(p)$. Then $\rk(V)\equiv 0\mod(p)$.

\vspace{2mm}

\noindent
(iv) Let $E$ be an exceptional object in $D(G(2,2^r))$. 
Then $\rk(E)\equiv m(V)+1\mod(2)$. 
\end{thm}

\Pf .
(i) By Proposition \ref{central-orth-prop}, we have an orthogonal decomposition 
$$K_{t_0}=\bigoplus_{m\in\Z/n\Z}K_{t_0}(m),$$
where $K_{t_0}=K_0^T(G(k,n))\otimes_R\Z[\sqrt[n]{1}]$,
and $K_{t_0}(m)$ is the $\Z[\sqrt[n]{1}]$-span of the classes of 
$E_i$ such that $m(E_i)\equiv m\mod(n)$.
Tensoring with $\OO(1)$ gives an isomorphism $K_{t_0}(m)\to K_{t_0}(m+k)$ of 
$\Z[\sqrt[n]{1}]$-modules.
Since $k$ is relatively prime to $n$, this implies that each 
$K_{t_0}(m)$ has rank $s/n$ over $\Z[\sqrt[n]{1}]$.

\noindent
(ii) By Corollary \ref{rk-cor},
it is enough to check that every bundle from the exceptional
collection $(\Sigma^{\la}\UU)$ described above has rank prime to $p$.
To this end we can use the hook-content formula to check that
the dimension of the irreducible representation of $\GL_k$ associated with partitions
$p-k\ge\la_1\ge\ldots\ge\la_k\ge 0$ is not divisible by $p$. Indeed, this dimension equals
$$s_{\la}=\prod_{x\in\la}\frac{k+c(x)}{h(x)},$$
where for every point $x=(i,j)$ of the Young diagram of $\la$, $h(x)$ is the hook length of $x$ and
$c(x)=j-i$ is the content of $x$ (see \cite{MacD},I.3, Ex.4). But our conditions on $\la$ imply that
$h(x)\le p-1$ and $k-i\le c(x)+k\le p-i$ for $x=(i,j)$, so all these numbers are relatively prime to $p$.
The last assertion follows from part (i) together with \eqref{m-deg-eq}.

\noindent
(iii) This follows immediately from \eqref{m-deg-eq}.

\noindent
(iv) By Corollary \ref{rk-cor}, 
it is enough to check this for bundles from the exceptional collection $(\Sigma^{\la}\UU)$ which
in this case consists of the symmetric powers of $\UU$ tensored with some line bundles.
\ed


\bigskip

{\it Quadrics}

\medskip

Let $Q^n$ denote the smooth quadric of dimension $n$, where $n\ge 3$.
It is a homogeneous space of the form $G/P$, where $G=\Spin(n+2)$. Recall
that we have a surjective homomorphism $G=\Spin(n+2)\to\SO(n+2)$ with the kernel of order $2$.
Let us denote by $Z_0\sub G$ this kernel. 

First, let us consider the case when $n$ is even. 
Then the group $G$ is simply connected of type $D_k$, where $n+2=2k$ ($k\ge 3$), 
$P$ is the maximal parabolic associated with the root $\a_1$. 
As in Proposition \ref{cases-prop}, we identify the standard maximal torus $T\sub G$ with the group 
$$\{(x_1,\ldots,x_k;x)\in\G_m^{k+1}\ |\ x^2=\prod x_i\}$$
in such a way that the projection to the coordinate $x_i$ corresponds to the character $\vareps_i$,
while the projection to the coordinate $x$ corresponds to the character $(\sum \vareps_i)/2$ of $T$.
Under this identification the projection to the first $k$ coordinates is exactly the map from $T$
to the maximal torus in $\SO(2k)$. Hence, $Z_0$ is generated by the element $z_0\in T$
that has all $x_i=1$ and $x=-1$. The Weyl group action on the weight lattice can permute $\vareps_i$'s 
and can multiply an even number of them by $-1$. Hence,
its action on $T$ is generated by permutations of coordinates $x_i$ together with the operator
$$(x_1,\ldots,x_k,x)\mapsto (x_1^{-1},x_2^{-1},x_3,\ldots,x_k,xx_1^{-1}x_2^{-1}).$$

We denote by $\OO(1)$ the $G$-equivariant line bundle on $G/P$ associated with the character
$\vareps_1$ of $T$ (it is an ample generator of the Picard group).
Note that $Z_0$ acts trivially on the line bundle $\OO(1)$.

We will use a slightly different element $t_0$ than in Proposition \ref{cases-prop}. Namely,
we set
$$t_0=(1,\zeta_n,\ldots,\zeta_n^{k-1};x_0),$$ 
where $x_0$ is a square root of $\zeta_{n}^{k(k-1)/2}$ 
(recall that $n=2k-2$). The assumption $(\star)$ is still satisfied 
for this element.

Let us denote by $w_1\in W$ the element such that
$$w_1(x_1,\ldots,x_k;x)=(x_1^{-1},x_2,\ldots,x_{k-1},x_k^{-1};xx_1^{-1}x_k^{-1}).$$
Then we have
\begin{equation}\label{quadric-w1-eq}
w_1(t_0)=z_0t_0.
\end{equation}

On the other hand, consider the element $w_2\in W$ such that
$$w_2(x_1,\ldots,x_k;x)=(x_k,x_{k-1}^{-1},\ldots,x_2^{-1},x_1^{\eps};x^{-1}x_1^{\frac{\eps+1}{2}}x_k),$$
where $\eps=(-1)^k$. Then we have
\begin{equation}\label{quadric-w2-eq}
w_2(t_0)=zt_0,
\end{equation}
where $z\in T$ is an element in the center of $\Spin(2k)$:
$$z=(-1,\ldots,-1:-\zeta_n^{-k(k-1)/2}).$$
Note that $z^2=z_0^k$.

Let $N$ be the order of $t_0$ (where $N|2n)$). Then
$t_0$ is a $\Q(\sqrt[N]{1})$-point of the torus $T$. Even though the set
$S=\{\zeta_n,\ldots,\zeta_n^{k-2}\}$ is not invariant under the action of $\Gal(\Q(\sqrt[N]{1})/\Q)$,
we claim that the conditions of Proposition \ref{Galois-prop} are still satisfied for $t_0$.
Namely, for every $\si\in \Gal(\Q(\sqrt[N]{1})/\Q)$ and every $s\in S$ we have either
$\si(s)\in S$ or $\si(s)^{-1}\in S$. Thus, we have a well defined set of signs $s_i=\pm 1$,
$i=2,\ldots,k-1$, and a permutation $w$ of $2,\ldots,k-1$, such that 
$$\si(t_0)=(1,\zeta_n^{s_2 w(2)},\ldots,\zeta_n^{s_{k-1}w(k-1)},-1;\si(x_0)).$$
Thus, there exists a unique element $\wt{w}_{\si}\in W$ such that $\wt{w}_{\si}(\vareps_k)=\vareps_k$ and
$$\si(t_0)=\wt{w}_{\si}(t_0)z_0^{\pi(\si)}$$
for some $\pi(\si)\in\Z/2\Z$. It follows that 
$\pi:\Gal(\Q(\sqrt[N]{1}/\Q)\to\Z/2\Z$ and $\si\mapsto \wt{w}_{\si}$
are group homomorphisms.
Now we set $w_{\si}=\wt{w}_{\si} w_1^{\pi(\si)}$. Since the elements
$\wt{w}_{\si}$ commute with $w_1$, we get that the map $\si\mapsto w_{\si}$ a group
homomorphism. From the above equation and from \eqref{quadric-w1-eq} we get
$$\si(t_0)=w_{\si}(t_0).$$

Recall (see \cite{Kap2}) that we have a full exceptional collection on $Q^{2k-2}=G/P$ consisting
of the $2k$ bundles $(\OO, S_+, S_-, \OO(1), \ldots, \OO(2k-3))$, where $S^{\pm}$ are the spinor bundles. The center $Z_G$ of $G=\Spin(2k)$ contains a subgroup $Z_0$ (the kernel of the homomorphism to $\SO(2k)$), and $Z_G/Z_0\simeq\Z/2\Z$. Let $\chi_0:Z_G\to\{\pm 1\}$ denote the unique nontrivial character of $Z_G$, that has trivial restriction to $Z_0$. 
Note that $Z_G$ acts on $\OO(1)$ through $\chi_0$.
Let $\chi_{\pm}:Z_G\to\G_m$ denote the characters with which the center acts on the spinor
bundles $S_{\pm}$. These characters are nontrivial on $Z_0$ and we have $\chi_+=\chi_0\chi_-$. 
The characters $\chi_0$ and $\chi_{\pm}$ are all nontrivial characters of $Z_G$.
Now Propositions \ref{central-orth-prop} and \ref{Galois-prop} 
give the following result.


\begin{cor}\label{even-quadric-cor}  
Let $V$ be a central object in $D(Q^{2k-2}/\Spin(2k))$, where $k\ge 3$. 

\noindent(i) If $V$ has trivial central character (resp., central character $\chi_0$) then the class 
$v(V)\in K_0^T(Q^{2k-2})\otimes_R\Z[\sqrt[N]{1}]$
is a linear combination with integer coefficients of the classes 
$v(\OO),v(\OO(2)),\ldots, v(\OO(2k-4))$ (resp., $v(\OO(1)),v(\OO(3)),\ldots,v(\OO(2k-3))$).

\noindent
(ii) If the central character of $V$ is $\chi_+$ (resp., $\chi_-$) then $v(V)$ is an integer multiple
of $v(S_+)$ (resp., $v(S_-)$).

\noindent
(iii) Assume now that $2k-2=2^r$ and the central character of $V$ is nontrivial. Then 
$\chi(V)\equiv 0\mod(2)$. If the central character of $V$ is $\chi_{\pm}$ then $\rk(V)\equiv 0\mod(2)$.
\end{cor}

Note that the last assertion in (iii) follows from (ii) and the fact that the rank of the spinor 
bundles $S_{\pm}$ is $2^{k-2}$, which is even since $k\ge 3$.

As in the case of Grassmannians, using the concrete full exceptional collection on $Q^{2k}$
we derive some information about arbitrary exceptional collections.

\begin{prop}\label{even-quadric-prop}
(i) Let $(E_1,\ldots, E_{2k})$ be a full exceptional collection in $D(Q^{2k-2})$, where $k\ge 3$. 
Then exactly $k-1$ objects from the collection have trivial central character,
$k-1$ objects have central character $\chi_0$, and the remaining two have central characters
$\chi_+$ and $\chi_-$.

\noindent
(ii) Let $E$ be an exceptional object in $D(Q^{2^r})$, where $r\ge 2$. Let us equip $E$ with
a $\Spin(2k)$-equivariant structure, where $k=2^{r-1}+1$.
Then $\rk E$ is odd iff the action of $Z_0$ on $E$ is trivial.
\end{prop}

\Pf . (i) Using equations \eqref{quadric-w1-eq},\eqref{quadric-w2-eq} and Proposition 
\ref{central-orth-prop}, we see that $K_0^T(Q^{2n})\otimes_R \Q[\sqrt[N]{1}]$ has a decomposition
into $4$ summands corresponding to different characters of $Z_G$. Now the assertion
follows from the form of the full exceptional collection on $Q^{2k}$.

\noindent
(ii) By Corollary \ref{rk-cor}, it is enough to check that this is true for the bundles of our exceptional
collection. 
\ed

Now let us consider odd-dimensional quadrics. These are homogeneous spaces of the form
$G/P$, where $G=\Spin(2k+1)$ is the simply connected group of type $B_k$ ($k\ge 2$),
$P$ is the maximal parabolic associated with the root $\a_1$ (the dimension of the quadric
equals $2k-1$). The maximal torus $T\sub G$ has the same description as in the case of type
$D_k$. However, the Weyl group is now bigger: it is generated by the group $S_k$
permuting the coordinates $x_i$'s
together with the involution 
\begin{equation}\label{w1-odd-eq}
(x_1,\ldots,x_k;x)\mapsto (x_1^{-1},x_2,\ldots,x_k;xx_1^{-1}).
\end{equation}

By Proposition \ref{cases-prop}, the following element satisfies the condition $(\star)$:
$$t_0=(-1,-\zeta_{2k-1},-\zeta_{2k-1}^2,\ldots,-\zeta_{2k-1}^{k-1};x_0),$$
where $x_0=i^k\zeta_{2k-1}^{k^2(k-1)/2}$.

The center of $\Spin(2k+1)$ coincides with $Z_0=\ker(\Spin(2k+1)\to\SO(2k+1))$. Its only nontrivial
element is $z_0=(1,\ldots,1;-1)\in T$.
Let $w_1\in W$ be the element acting on $T$ by the involution \eqref{w1-odd-eq}. Then
$$w_1(t_0)=z_0t_0.$$

As in the case of even-dimensional quadrics, we can check that
the conditions of Proposition \ref{Galois-prop} are satisfied for $t_0$.
Namely, let $N$ be the order of $t_0$. Then for every $\si\in \Gal(\Q(\sqrt[N]{1})/\Q)$
we have a well defined set of signs $s_i=\pm 1$,
$i=2,\ldots,k$, and a permutation $w$ of $2,\ldots,k$, such that 
$$\si(t_0)=(-1,\zeta_{2k-1}^{s_2 w(2)},\ldots,\zeta_{2k-1}^{s_{k}w(k)};\si(x_0)).$$
Thus, as before we have group homomorphisms
$\pi:\Gal(\Q(\sqrt[N]{1})/\Q)\to\Z/2\Z$ and $\si\mapsto \wt{w}_{\si}$, where
$\wt{w}_{\si}(\vareps_1)=\vareps_1$ and
$$\si(t_0)=\wt{w}_{\si}(t_0)z_0^{\pi(\si)}$$
Hence, setting $w_{\si}=\wt{w}_{\si} w_1^{\pi(\si)}$ we get
$$\si(t_0)=w_{\si}(t_0).$$

We have a full exceptional collection on $Q^{2k-1}=G/P$ consisting
of the $2k$ bundles $(\OO, S, \OO(1), \ldots, \OO(2k-2))$, where $S$ is the spinor bundle
(see \cite{Kap2}). The center 
of $\Spin(2k+1)$ acts trivially on line bundles, and nontrivially on $S$.
Thus, Propositions \ref{central-orth-prop} and \ref{Galois-prop} 
give the following result.

\begin{cor}\label{odd-quadric-cor}  
Let $V$ be a central object in $D(Q^{2k-1}/\Spin(2k+1))$. 

\noindent(i) If $V$ has trivial central character then the class 
$v(V)\in K_0^T(Q^{2k-1})\otimes_R\Z[\sqrt[N]{1}]$
is a linear combination with integer coefficients of the classes 
$v(\OO),v(\OO(1)),\ldots, v(\OO(2k-2))$.

\noindent
(ii) If $V$ has a nontrivial central character then $v(V)$ is an integer multiple
of $v(S)$.
\end{cor}

As in the case of even-dimensional quadrics, we also obtain some information on central
characters of objects in an arbitrary full exceptional collection.

\begin{prop}\label{odd-quadric-prop}
Let $(E_1,\ldots, E_{2k})$ be a full exceptional collection in $D(Q^{2k-1})$. 
Then one of the objects $E_i$ has a nontrivial central character, while the remaining
$2k-1$ have trivial central character.
\end{prop}

\bigskip

{\it Maximal isotropic Grassmannians}

\medskip

Let us denote by $OG(k,n)$ (resp., $SG(k,n)$)
the (connected component of) the Grassmannian of isotropic
$k$-dimensional subspaces in an $n$-dimensional orthogonal (resp., symplectic) vector space.
The maximal isotropic Grassmannians are $OG(k,2k)$, $OG(k,2k+1)$ and $SG(k,2k)$.
Note that there is an isomorphism $OG(k,2k+1)\simeq OG(k+1,2k+2)$.
At present, the existence of full exceptional collections is not known in these cases (except
in small dimensions). However, it is known that $K^T_0(X)$ is a free $R(T)$-module (in fact, this is always true for generalized flag varieties $X=G/P$, see \cite{KK}).

Let us first consider the orthogonal case: $X=OG(k,2k)$. Then $X=G/P$, where $G=\Spin(2k)$,
$P$ is the maximal parabolic subgroup associated with $\a_k$. 
Let us consider
the same element $t_0=(1,\zeta_{2k-2},\ldots,\zeta_{2k-2}^{k-1};x_0)$ as for the even-dimensional
quadric. By Proposition \ref{cases-prop}, it satisfies the condition $(\star)$.
Hence, Propositions \ref{central-orth-prop} and \ref{Galois-prop} are still applicable in this case.
Note that the generator of the Picard group $\OO(1)$ in this case corresponds to the character
$\om_k=(\sum \vareps_i)/2$, so $Z_0$ acts on it nontrivially. Thus, from Proposition
\ref{central-orth-prop} we derive the following result.

\begin{cor}
(i) Let $N$ be the order of $t_0$ (note that $N|4(k-1)$). Then we have an orthogonal decomposition
$$K_0^T(OG(k,2k))\otimes_R \Z[\sqrt[N]{1}]=M\oplus M\cdot[\OO(1)],$$
where $M\sub K_0^T(OG(k,2k))\otimes_R \Z[\sqrt[N]{1}]$ is the $\Z[\sqrt[N]{1}]$-span
of the classes of $\SO(2k)$-equivariant bundles. If in addition $k$ is odd then
there is an orthogonal decomposition
$$M=M_0\oplus M_0\cdot[\OO(2)],$$
where $M_0\sub M$ is the $\Z[\sqrt[N]{1}]$-span of the classes of $\SO(2k)/\{\pm 1\}$-equivariant
bundles.

\noindent
(ii) Let $V$ be a central object in $D(OG(k,2k)/\Spin(2k))$ 
with a nontrivial central character.
Assume that $k=2^r+1$. Then $\chi(V)\equiv 0\mod(2)$.
\end{cor}

In the symplectic case we have $X=SG(k,2k)=\Sp(2k)/P$, where $P$ is the maximal parabolic
associated with $\a_k$. As in Proposition \ref{cases-prop}, let us consider the element
$$t_0=(\zeta_{2k+2},\zeta_{2k+2}^2,\ldots,\zeta_{2k+2}^k)$$
satisfying the condition $(\star)$.
Let also $w_1$ be the element of the Weyl group sending
$(x_1,x_2,\ldots,x_k)$ to $(x_k^{-1},\ldots,x_2^{-1},x_1^{-1})$. Then we have
$$w_1(t_0)=zt_0,$$
where $z=(-1,\ldots,-1)\in T$ is the generator of the center of $\Sp(2k)$.
The generator of the Picard group $\OO(1)$ corresponds to the character $\sum \vareps_i$.
Hence, the element $z$ acts on $\OO(1)$ as $(-1)^k$.
As before, Proposition \ref{central-orth-prop} gives the following result.

\begin{cor}(i) Let $N=2k+2$. We have an orthogonal decomposition
$$K_0^T(SG(k,2k))\otimes_R \Z[\sqrt[N]{1}]=M_+\oplus M_{-},$$
where $M_+$ (resp., $M_{-}$) is the $\Z[\sqrt[N]{1}]$-span of the classes of
$\Sp(2k)$-equivariant bundles with trivial (resp., nontrivial) central character.
If $k$ is odd then $M_{-}=M_+\cdot[\OO(1)]$.

\noindent
(ii) Let $V$ be a central object in $D(SG(k,2k)/\Sp(2k))$ with a nontrivial
central character. Assume that $k=2^r-1$. Then $\chi(V)\equiv 0\mod(2)$.
\end{cor}

The first interesting examples of maximal isotropic Grassmannians are $SG(3,6)$
and $OG(5,10)$ (note that $SG(2,4)$ and $OG(4,8)$ are smooth
quadrics). There exists a full exceptional collection on
$SG(3,6)$ (resp., $OG(5,10)$) consisting of line bundles and vector bundles
of rank $3$ (resp., $5$), see \cite{Kuznetsov-hyp}, 6.2 and 6.3 (also \cite{Sam} in the case of $SG(3,6)$).
Hence, Corollary \ref{rk-cor} leads to the following result.

\begin{thm}\label{isotropic-thm} 
Every exceptional object on $SG(3,6)$ (resp., $OG(5,10)$) has an odd rank.
\end{thm} 

\bigskip

{\it Hirzebruch surfaces and products}

\medskip

For an integer $n$ let us consider the ruled surface 
$$F=F_n=\P(\OO_{\P^1}\oplus\OO_{\P^1}(n))\stackrel{\pi}{\to}\P^1.$$
We have a natural $\GL_2$-action on $F$ induced by its action on $\P^1$.
In addition, let us equip it with the fiberwise action of $\G_m$ that acts trivially on 
$\OO_{\P^1}$ and by the identity character on $\OO_{\P^1}(n)$.
The maximal torus $T=\G_m^2\sub\GL_2$ 
acts on $F$ with $4$ stable points: two on the fiber $\pi^{-1}(p_1)$
and two on $\pi^{-1}(p_2)$ (where $p_1$ and $p_2$ are $T$-stable points on $\P^1$). 
The tangent bundle $\TT_F$ to $F$ fits into the exact sequence
$$0\to \TT_{\pi}\to \TT_F\to\pi^*\TT_{\P^1}\to 0,$$
where $\TT_{\pi}\simeq \OO_F(2)\otimes\pi^*\OO_{\P^1}(n)$.
Hence, the weights of the action of $(x_1,x_2,u)\in T\times \G_m$ 
on the tangent spaces to the $T$-stable points are:
(i) for two points in $\pi^{-1}(p_1)$: $(x_1^nu,x_1/x_2)$ and $(x_1^{-n}u^{-1},x_1/x_2)$;
(ii) for two points in $\pi^{-1}(p_2)$: $(x_2^nu,x_2/x_1)$ and $(x_2^{-n}u^{-1},x_2/x_1)$.
Thus, in the case when $n$ is even
the condition $(\star)$ of Theorem \ref{main-thm} is satisfied for an element 
$$t_0=
\cases
(i,i^{n+1},1) &\text{ if }n\equiv 2(4),\\
(i,i^{n-1},-1) &\text{ if }n\equiv 0(4)\endcases
$$
of order $4$. Thus, we deduce from Corollary \ref{rk-cor} the following result.

\begin{cor}\label{odd-prod-cor} Let 
$X=F_{2n_1}\times\ldots\times F_{2n_r}\times \P^{2^{k_1}-1}\times\ldots\times\P^{2^{k_s}-1}$.
Then for every exceptional object $E$ in $D(X)$ the class of $E$ in $K_0(X)\otimes\Z/2\Z$
coincides with the class of one of the line bundles. In particular, $\rk(E)$ is odd.
\end{cor}

\begin{rem}
Let $X$ be a smooth projective variety that admits a full exceptional collection $(E_i)$
(where $E_i\in D(X)$) consisting of objects of odd rank. Assume also that
$$\chi(x,y)\equiv\chi(y,x)\mod 2$$
for all $x,y\in K_0(X)$.
Then expressing $[V]\in K_0(X)$ in terms of the basis $([E_i])$ one immediately
checks that
$$\chi(V,V)\equiv\rk(V)\mod(2).$$
In particular, every exceptional object on $X$ has odd rank.
The class of varieties with above properties is closed under products and includes
projective spaces of dimension $2^k-1$, Hirzebruch surfaces $F_n$ for even $n$,
and varieties $OG(5,10)$ and $SG(3,6)$.
\end{rem}

\bigskip

{\it Relative projective spaces}

\medskip

Let $S$ be a smooth projective variety.
We consider the relative projective space $\P^{n-1}_S=\P^{n-1}\times S$.
Recall that $p_1,\ldots,p_{n}\in\P^{n-1}$ denote the $T$-stable points, where 
$T\sub\GL_{n}$ is the standard maximal torus, and the action of $T\simeq\G_m^n$
on the fiber of $\OO(1)$ at $p_i$ is given by the projection to the $i$-th coordinate.
Consider the element $t_0\in T$ given by \eqref{GL-t0-eq} and the corresponding
homomorphism $R=R(T)\to\Z[\sqrt[n]{1}]$.
For an object $V\in D(\P^{n-1}_S/T)$ let us denote by $v(V)$ the corresponding
class in $K_0^T(\P^{n-1}_S)\otimes_{R}\Z[\sqrt[n]{1}]$. 

As before, for a central object in $D(\P^{n-1}_S/\GL_n)$ with 
the central character $\la\mapsto\la^m$ we set $m(V)=m\mod(n)$. 
By Lemma \ref{equiv-lem}, every exceptional object $E$ of $D(\P^{n-1}_S)$
has a $\GL_n$-equivariant structure, unique up to tensoring with a character. 
Hence, $m(E)\in\Z/n\Z$ is well defined.

\begin{thm}\label{rel-proj-thm} 
(i) For a central object $V\in D(\P^{n-1}_S/\GL_n)$ set
$$\tau(V):=\tr(t_0,V|_{p_1\times S})\in K_0(S)\otimes \Z[\sqrt[n]{1}].$$
Then $\tau(V)$ belongs to $K_0(S)\sub K_0(S)\otimes \Z[\sqrt[n]{1}]$, and 
\begin{equation}\label{epsV-eq}
v(V)=\tau(V)\cdot v(\OO(m(V)))
\end{equation}
in $K_0^T(\P^{n-1}_S)\otimes_{R}\Z[\sqrt[n]{1}]$.

\noindent
(ii) Let $E$ be an exceptional object of $D(\P^{n-1}_S)$. Let us equip $E$ with a $\GL_n$-equivariant
structure and consider the element $\tau(E)\in K_0(S)$. Then
$$\chi_S(\tau(E),\tau(E))=1,$$
where $\chi_S$ denotes the Euler form on $K_0(S)$. 
If $(E_1,E_2)$ is an exceptional pair then either $m(E_1)\not\equiv m(E_2)$ 
or $\chi_S(\tau(E_2),\tau(E_1))=0$.
If $(E_1,\ldots,E_N)$ is a full exceptional collection in $D(\P^{n-1}_S)$,
then for each $m\in\Z/n\Z$ exactly $N/n$ objects $E_i$ have $m(E_i)\equiv m$,
and the corresponding $N/n$ elements $\tau(E_i)$ form a semiorthogonal basis of $K_0(S)$.

\noindent
(iii) Assume now that $n=p^r$, where $p$ is prime. Then for a central object $V\in D(\P^{n-1}_S/\GL_n)$
one has the following congruence in $K_0(\P^{n-1}_S)\otimes\Z/p\Z$:
$$[V]\equiv [V|_{p_1\times S}]\cdot [\OO(m(V))].$$
If $E$ is an exceptional object of $D(\P^{n-1}_S)$ then 
$$\chi_S(E|_{p_1\times S},E|_{p_1\times S})\equiv 1\mod(p).$$
If $(E_1,E_2)$ is an exceptional pair with $m(E_1)\not\equiv m(E_2)\mod(n)$ then
$\chi(E_1,E_2)\equiv 0\mod(p)$. In the case $m(E_1)\equiv m(E_2)\mod(n)$ we have
$\chi_S(E_2|_{p_1\times S},E_1|_{p_1\times S})\equiv 0\mod(p)$.
\end{thm}

\Pf . 
(i) By the projective bundle theorem, we have a decomposition
$$K_0^T(\P^{n-1}_S)=\bigoplus_{i=0}^{n-1} [\OO(i)]\cdot K_0(S)\otimes R.$$
Hence, 
$$K_0^T(\P^{n-1}_S)\otimes_{R}\Z[\sqrt[n]{1}]=\bigoplus_{i=0}^{n-1}M_i,$$ 
where $M_i$ is the $\Z[\sqrt[n]{1}]$-submodule generated by $K_0(S)\cdot v(\OO(i))$.
Applying Proposition \ref{central-orth-prop} we see that $M_i$ coincides with the
submodule generated by the classes of $\GL_n$-equivariant bundles $V$ with $m(V)\equiv i\mod(n)$. In particular, the submodules $M_i$ and $M_j$ are
orthogonal with respect to the form $\Tr(t_0,\chi^T(\cdot,\cdot))$ 
for $i\neq j$. Furthermore, the restriction to $p_1\times S$ gives the inverse of
the natural isomorphism
$$K_0(S)\otimes\Z[\sqrt[n]{1}]\wt{\to} M_i:x\mapsto x\cdot[\OO(i)]$$
(since $t_0$ acts trivially on the fiber of
$\OO(1)$ at $p_1$). This immediately implies \eqref{epsV-eq}.
Note that $\tau(V)$ is the specialization at $t_0$ of the class of $V|_{p_1\times S}$ in
$K_0^T(S)=K_0(S)\otimes R$. Since $V$ is $\GL_n$-equivariant, the class 
$V|_{p_1\times S}$ is preserved under the action of $S_{n-1}\sub S_n$ that stabilizes $p_1$.
Therefore, its specialization at $t_0$ is invariant under the Galois action on $\Z[\sqrt[n]{1}]$
(permuting the nontrivial $n$-th roots of unity), so $\tau(V)\in K_0(S)$.

\noindent
(ii) If $E$ is exceptional then we have
$1=\Tr(t_0,\chi^T(E,E))$. Hence, using \eqref{epsV-eq} we immediately deduce that
$\chi_S(\tau(E),\tau(E))=1$. In the same way we get the vanishing of $\chi_S(\tau(E_2),\tau(E_1))$
in the case when $(E_1,E_2)$ is an exceptional pair with $m(E_1)\equiv m(E_2)$. This easily implies
the last assertion.

\noindent
(iii) This follows from (i) and (ii) using the homomorphism $\rho:\Z[\sqrt[n]{1}]\to\Z/p\Z$.
\ed



\end{document}